\documentclass[12pt]{amsart}
\usepackage{amssymb}

\newcommand\datver[1]{\def\datverp%
 {\par\boxed{\boxed{\text{Version: #1; Run: \today}}}}}
\datver{1.1X; Revised: 03-11-98}

\newcommand\B{\mathcal B}

\newcommand\C{\mathbb C}
\newcommand\R{\mathbb R}
\newcommand\HT{{\mathcal H}_t}

\newcommand\DD{\mathbb D}

\def \B* {\overline {B}^*}

\newcommand\Id{\operatorname{Id}}

\newtheorem{theorem}{Theorem}

\newtheorem{proposition}{Proposition}
\newtheorem{corollary}{Corollary}

\theoremstyle{definition}
\newtheorem{definition}{Definition}

\newtheorem{remark}{Remark}
\begin{document}

\title[Index of ${\Gamma}$-equivariant Toeplitz operators]{Index of
  ${\Gamma}$-equivariant Toeplitz operators}

\author[R. Nest]{Ryszard Nest}
    \address{Copenhagen}
       \email{rnest@math.ku.dk}
\author[F. Radulescu]{Florin Radulescu}
    \address{ University of Iowa, Math. Dept.,Iowa City
Iowa 52242, USA}


\begin{abstract}

 Let  $\Gamma$ be a discrete subgroup
of $PSL(2,\Bbb R)$ of infinite covolume with infinite
conjugacy classes. Let ${\HT}$ be the
Hilbert space consisting of analytic functions
 in $L^2(\DD,(\text{Im\ } z)^{t-2}\text{d}\overline{z}\text{d}z)$ 
 and let, for $t>1$, $\pi_t$ denote the
 corresponding projective unitary representation of $PSL(2,\Bbb R)$ on this
Hilbert space. We denote by $\mathcal A_t$ the $II_\infty$ factor
given by the commutant of 
$\pi_t(\Gamma)$ in $B({\HT})$.  Let $F$ denote a fundamental
domain for $\Gamma$ in $\DD$ and assume that $t >5$. $\partial
M=\partial \DD \cap \overline{F}$ is given the topology of disjoint
union of its connected components.

Suppose that $f$ is a continuous $\Gamma$-invariant function on $\DD$
whose restriction to $F$
extends to a continuous function on  $\overline{F}$ and such that
$f|_{\partial M}$ is an invertible element of $C_0
(\partial M)\mbox{\~{}}$. 
 Let $T^t_f={\mathcal P}_t  {\mathcal M}_{{f}} {\mathcal P}_t 
$ denote the Toeplitz operator with symbol $f$.
Then $T^t_f$ is Fredholm, in the Breuer sense, with respect to the 
$II_\infty $ factor $\mathcal A_t$ and, moreover, its Breuer index is
equal to the total winding number of $f$ on $\partial M$. 

\end{abstract}

\maketitle

\tableofcontents

\section{Introduction}

In this paper we study equivariant Toeplitz operators acting on the
Hilbert space 
${\HT}$ consisting of all square summable analytic functions in 
$$
L^2(\DD,(\text{Im\ } z)^{t-2}\text{d}\overline{z}\text{d}z).
$$ 
Let us first recall that the classical theory of Toeplitz operators in
the unit disc yields an extension of C*-algebras
$$
0\rightarrow {\mathcal K} \rightarrow {\mathcal T} \rightarrow C(\partial
\DD ) \rightarrow 0,
$$
where ${\mathcal K}$ denotes the algebra of compact operators on
${\mathcal H}_2$ and $\mathcal T$ the Toeplitz C*-algebra generated by
compressions $T_f$ to ${\mathcal H}_2$ of multiplication operators (by
$f$'s fom $C(\overline{\DD })$). In particular, for $f|_{\partial \DD}$ invertible, 
$T_f$ is Fredholm and the boundary map for the K-theory six
term exact sequence of this extension is equivalent, 
via index theorem for Toeplitz operators, to the equality:  
$$
\mbox{Index } (T_f )
=\mbox{winding number of $f|_{\partial \DD}$}.
$$
As it turns out, all of these facts admit suitable generalisation to
the equivariant case.

\medskip

Let
 $\Gamma$ be a fuchsian subgroup of $PSL(2,\mathbb R)$, which has infinite 
conjugacy
 classes and is of infinite covolume. Recall that the action of
 PSL$(2,\mathbb R)$ on $\DD$ by fractional linear transformations
 lifts to projective
 unitary representations  of 
PSL$(2,\mathbb R)$ on these Hilbert spaces (cf. (\cite {Sa}, 
 \cite {Pu})) and 
 the commutant of $\pi_t(\Gamma)$
is a $II_\infty$ factor. We will denote by  $\mathcal A_t$ the
commutant $\pi_t (\Gamma )^{'} \cap {\mathcal B}(\HT )$ and by $\tau$
the normal positive non-zero trace on $\mathcal A_t$.  
 If $\sigma_t$ denotes 
 the 2-group cocycle  
corresponding to the projective 
unitary representation $\pi_t$, then $\mathcal A_t$
is isomorphic  (\cite{Ra}) to $\mathcal L(\Gamma,\sigma_t)\otimes B(K)$,
where $K$ is an infinite dimensional separable Hilbert space and 
$\mathcal L(\Gamma,\sigma_t)$ is the twisted group von Neumann algebra
of $\Gamma$.

Let $F$ denote a fundamental domain for the action of $\Gamma$ on
 $\DD$. We will denote by $M$ the quotient space $\DD /\Gamma$ and by
 $\partial M$ its boundary: 
$$
\partial M =  (\partial \DD \cap \overline {F} )/\Gamma
$$
equipped with the topology of the disjoint union of its connected
components, i. e. of a countable union of disjoint circles. In
particular $\overline{M} = \overline{F}/\Gamma = M \cup \partial M$
inherits a structure of locally compact space. $C_0
(\partial M)\mbox{\~{}}$ denotes the unitalisation of the C*-algebra
of continuous, vanishing at infinity functions on $\partial M$.

For any  continuous $\Gamma$-invariant function  $f$ on $\DD$ which
  extends to a continuous function on $\overline{M}$ we denote by $T_f$ the
Toeplitz operator 
on ${\HT}$ with symbol $f$, i. e. the compression to $\HT$ of the
  operator of multiplication by $f$ on $L^2(\DD,(\text{Im\ }
  z)^{t-2}\text{d}\overline{z}\text{d}z) $. Because of the $\Gamma$
  invariance of the 
symbol such a $T_f$ belongs to $\mathcal 
A_t=\{\pi_t(\Gamma)\}'$
 (\cite{Ra}).

Suppose that $t >5$ and
that $f|_{\partial M}$ is an invertible element of $C_0
(\partial M)\mbox{\~{}}$. We prove below that $T_f$ is Fredholm (in
the sense of Breuer (\cite {Bre}),
in $\mathcal A_t$. Moreover
the Breuer index is (in analogy with the classical case) equal to 
the winding number of $f|_{\partial M}: \partial M \rightarrow \C
\setminus \{ 0 \}$. 
This can be seen as an analogue of Atiyah's index formula  for
coverings (\cite{At}).

\medskip

The organisation of the paper is as follows.

\medskip

In Section 2 we gather some
more or less known results about nuclearity properties of Toeplitz
operators on $\HT$ and prove the main technical result:

\begin{em} 
Let $t >5 $ and $f,\ g\in L^{\infty}({\overline{\DD}})$ are given. Suppose
that $g\in C^{\infty} (\overline{\DD})$ and that 
$$
\mbox{inf}\ \{ ||z-\xi || | \mid z\in \mbox{supp $f$ and $\xi \in$supp $g$}
\} > \epsilon 
$$
for some positive number $\epsilon$. Then both $T_f T_g$ and $T_g T_f$
are of trace class and Tr$ ([ T_f ,T_g ]) =0$.
\end{em}

\noindent (cf. Theorem \ref{t2} ).

\medskip

In Section 3 we study the ${\mathcal L}^1 (\tau )$-properties of
commutators of Toeplitz operators with $\Gamma$-invariant symbol and
prove the following result.

\begin{em}
Suppose that $t>5$ and $f$ and $g$ are $\Gamma$-invariant functions on
$\DD$ which are smooth on the closure of a fundamental domain for
$\Gamma$.
Then both $[T_f ,T_g ]$ and $T_{fg}-T_f T_g$ are in ${\mathcal
M}\cap{\mathcal L}^1 (\tau)$ and
$$
\tau ([T_f ,T_g ] )=\frac{1}{2\pi i} \int_F df dg 
$$
\end{em}
(cf. Theorem \ref{t3} and the remarks following).

\medskip

Let ${\mathcal
T}_{\Gamma}$ be the C*-subalgebra of ${\mathcal A}_t$ generated by
Toeplitz operators $T_f$ with $f$ 
$\Gamma $-invariant and smooth on a fundamental domain for $\Gamma$,
${\mathcal K}_{\Gamma}$ be the C*-ideal generated by the ${\mathcal L}^1
(\tau )$-elements in ${\mathcal
T}_{\Gamma}$ and $M=\DD /\Gamma$. In Section 4 we construct the extension 
$$
\begin{array}{ccccccccc}
0&\rightarrow& {\mathcal K}_{\Gamma}& \rightarrow& {\mathcal
T}_{\Gamma}& \rightarrow &C(\partial M)& \rightarrow &0.\\
&&&&T_f&\rightarrow & f|_{\partial M}&&\\
\end{array}
$$
Let $\partial : K_1 (\partial M) \rightarrow K_0 ({\mathcal K}_{\Gamma} ) $ denote
the boundary map in K-theory associated to this extension. We prove that, 
for $T_f \in {\mathcal T}_{\Gamma}$ with 
$f$ invertible on the boundary of $M$, the following
equality holds:
$$
\langle \tau \ , \ \partial [T_f ] \rangle = \mbox{winding number of $f$ on
$\partial M$}
$$
(cf. Theorem 4).

\medskip

\begin{remark}
For notational simplicity we work throughout the paper with the case
when the number of boundary components of a fundamental  
domain of $\Gamma$ is finite. The only difference (except for
typografical complications) in the general case consists of replacing
the above extension of $C(\partial M)$ by ${\mathcal K}_{\Gamma}$ by
the extension: 
$$
0\rightarrow {\mathcal K}_{\Gamma} \rightarrow {\mathcal T}^0_{\Gamma}
\rightarrow 
C_0 (\partial M) \rightarrow 0
$$ 
where ${\mathcal T}^0_{\Gamma}$ stands for the (in general nonunital) C*-algebra
generated by $T_f$ with $f$ continuous on $\overline{F}$ 
and in with non-zero values on finitely many components of 
$\partial \DD \cap \overline{F}$.
\end{remark}

\medskip

The method of the proof are based on the equivariant Berezin's quantization
theory for such groups (\cite{Ra}, \cite{NN}). Let $F$ be, as above, a
fundamental domain 
for
the action of $\Gamma$ in $\DD$. Then for 
every $\Gamma-$ equivariant, bounded function $g$ on $\DD$, having compact 
support in the interior of $F$, 
the Toeplitz operator $T_g\in\mathcal A_t$ is in 
${\mathcal L}^1(\mathcal A_t)$ and has 
trace
equal to a universal constant times the integral
$\int_F g(z) (\text{Im\ })^{-2}\text{d}z\text{d}\overline z$

Moreover we will show that the commutator of two  Toeplitz opertors,
having symbols that are smooth  and continuous on the closure of $F$, 
belongs to trace ideal of the
$II_\infty$ factor. In particular if the symbol is invertible 
in the neighbourhood of the intersection
of the boundary of $\DD$ with the closure of $F$,
the operator is Fredholm in $\mathcal A_t$ in Breuer's sense 
 (\cite{Bre}).

To identify the Breuer index of such
 a Toeplitz operator we use the Carey-Pincus theory (\cite{CP2}.
Let us first recall the pertinent
facts. Given an operator $A\in {\mathcal A}_t$ such that 
$\tau [A^* ,A] < \infty$, the bilinear map
$$
\C [\overline{z},z]\ni P,Q \rightarrow \tau [P(A^* ,A),Q(A^* ,A))] 
$$ 
defines a cyclic one-cocycle on the algebra of polynomials (in two real
variables), of the form
$$
(P,Q)\rightarrow \int_{Z}\{ P,Q \} d\mu (\overline{z},z),
$$
where $Z$ is the of spectrum of the class of $A$ in 
$ {\mathcal A}_t /{\mathcal K}_{\Gamma}$.
$d\mu$, called the principal function of $A$, is a finite measure on
the complex plane
having the property that, for any connected component $Z^{'}$ 
of $\C \setminus Z$, 
$$
d\mu|_{Z^{'}} = - c \pi d\overline{z}dz
$$
with a constant $c$ equal to the value of $\tau$ on the index class of 
$\partial [A - \lambda] \in \mbox{K}_0 ( {\mathcal K}_\Gamma )$ 
for any $\lambda \in Z^{'}$. 
Hence one has to determine the pricipal function, 
given by the
$\tau$-values on commutators of polynomials in $T_f$
and its adjoint (like in the classical case in \cite{CP1}, \cite{HH}).

To deal with the computation of those we
apply the spatial theory
of von Neumann algebras (\cite{C-1}). One of its basic constructions 
gives an operator-valued weight
$E: {\mathcal B}(\HT)\rightarrow {\mathcal A}_t$ such that, 
for a trace-class operator $A_0$ in the domain of $E$, 
$\tau \circ E(A_0 )=Tr (A_0)$. 
This allowes one to replace the computation
of values of $\tau$ on $\Gamma$-invariant operators by  
computation of value of the classical trace on certain trace-class
operators on $\HT$. In fact we prove in Section 3 that commutators of 
$\Gamma$-equivariant Toeplitz operators are of the form $E(A_0)$ for 
$A_0$ given by a suitable (trace class) commutator of polynomials in
Toeplitz operators and hence the computation reduce to the classical
case.

But in this case the principal function for a pair of Toeplitz operators that
commute modulo the trace ideal in $B({\HT})$ is well understood - 
it is basically given by the fact that the index of $T_{z}$ is equal to one 
(\cite{Arazy}) and gives explicit formulas that lead to the results
stated above.


\section{Some results on Toeplitz operators}

Let $\DD$ denote the unit circle in the complex plane and set
\begin{equation}
d\mu_t (z)=(1-|z|^2)^t \frac{d\bar{z}dz}{(1-|z|^2)^2}
\end{equation}
We set
\begin{equation}\label{1}
{\mathcal H}_t =\{ f\in L^2 (\DD , d\mu_t )\mid \mbox{ $f$ holomorphic on }\DD 
\}
\end{equation}
As is well known, ${\mathcal H}_t$ is a closed subspace of $L^2 (\DD ,
d\mu_t )$
and the orthogonal projection ${\mathcal P}_t : L^2 (\DD , d\mu_t )\rightarrow
{\mathcal H}_t$
is called the Toeplitz projection. Given a function $f\in L^\infty (\DD ,
d\mu_t )$ we denote by ${\mathcal M}_f$ the operator of multiplication
by $f$ on $L^2 (\DD , d\mu_t )$ and set
\begin{itemize}
\item the Toeplitz operator associated to $f$: 
$$
T_f = {\mathcal P}_t {\mathcal M}_f {\mathcal P}_t 
$$

\item the Henkel operator associated to $f$:
$$
H_f = (1-{\mathcal P}_t ) {\mathcal M}_{\bar{f}} {\mathcal P}_t 
$$
\end{itemize}
Note, for future computations, that $T_f$ is an integral operator on
$L^2 (\DD , d\mu_t )$ with
integral kernel
\begin{equation}
K_f (z, \xi )= \frac{t-1}{2\pi i} \frac{f(\xi )}{(1-z\bar{\xi })^t}
\end{equation}
and
\begin{equation}
H_{\bar{f}}^* H_g =T_{fg} - T_f T_g.
\end{equation}

We will denote by $\delta$ the absolute value of the cosine of
hyperbolic distance on $\DD$, i. e.
\begin{equation}\label{2}
\delta (a,b)=\frac{(1-||a||^2)(1-||b||^2)}{|1-a\bar{b}|^2}.
\end{equation}

 As is well known,
\begin{equation}\label{3}
||T_f || = ||f||_\infty
\end{equation}
and
\begin{equation}\label{4}
||f||_{S_2}^2 =||H_f||^2_2 + ||H_{\bar{f}}||^2_2 = 
(\frac{t-1}{2\pi i})^2 \int_{\DD \times \DD}
|f(a)-f(b)|^2 \delta^t (a,b) d\mu_0 (a,b).
\end{equation}
In particular, since both $H_z$ and $H_{\bar{z}}$ are of finite rank,
the function $f(a,b)=(a-b)$ is square integrable with respect to the measure
$\delta^t (a,b) d\mu_0 (a,b)$ and hence all functions which are
Lipschitz with exponent one on $\DD$ have finite $S_2$-norm.

We let PSL(2,$\R $) act on $\DD$ by fractional linear transformations
and denote by $\pi_t$ the induced projective unitary representation on
$L^2 (\DD , d\mu_t)$ (and ${\mathcal H}_t$). Both $d\mu_0$ and $\delta
(a,b)$ are PSL(2,$\R $)-invariant, which gives a useful formula
\begin{equation}\label{5}
\int_\DD \delta^t (a,b) d\mu_0 (a) = \int_\DD \delta^t (a,0) d\mu_0 (a)
= \frac{4\pi }{t-1}.
\end{equation}

\bigskip 
The following result is probably well known to specialists, however,
since we do not
have a ready reference, so we will include the proof below.
\begin{theorem}\label{t1}
Let $f$ and $g$ belong to $C^{\infty}(\overline{\DD})$. Then $T_f T_g -
T_{fg}$ is a trace class operator and, moreover,
$$
Tr ([ T_f ,T_g ] )=(\frac{t-1}{2\pi i})^2 \int_{\DD \times
\DD}(f(a)g(b)-f(b)g(a))\delta^t (a,b) d\mu_0 (a,b)$$
$$=\frac{1}{2\pi i} \int_{\DD} dfdg =\frac{1}{2\pi i} \int_{\partial \DD} fdg.
$$
\end{theorem}
Proof. By the smoothness assumption, both $f$ and $g$ belong to the
$S_2$ class and hence $T_f T_g -
T_{fg}$ is a trace class operator. For the computation of the trace we
can just as well assume that both $f$ and $g$ are real-valued. To begin
with, for a real-valued function $f$, (\ref{4}) gives
$$
Tr (T_{f^2} -T_f^2 )=\frac{1}{2} (\frac{t-1}{2\pi i})^2 \int_{\DD \times
\DD} (f(a)-f(b))^2
\delta^t (a,b) d\mu_0
$$
and hence, by an application of the polarisation identity, 
$$
Tr (T_{fg}-T_{f} T_g )= (\frac{t-1}{2\pi i})^2  \int_{\DD \times \DD}
f(a)(g(a)-g(b))\delta^t (a,b) d\mu_0
$$
which implies immediately the first equality.

To get the second equality recall that, for a pair of (non-commutative)
polynomials $P(\bar{z},z)$ and $Q(\bar{z},z)$, the Carey-Pincus formula
holds:
$$
Tr ([ P(T^*_z ,T_z ),Q(T^*_z ,T_z )])=\frac{1}{2\pi i} \int_{\DD} dPdQ 
=\frac{1}{2\pi i} \int_{\partial \DD} PdQ,
$$
see \cite{CP1}. Since $P(T^*_z ,T_z )=T_P \ \mbox{mod}\ {\mathcal L}^1
(Tr) $, this implies the second equality for $f$ and $g$ polynomial.
Approximating arbitrary pair of smooth functions uniformly with their first
derivatives on $\overline{\DD}$ completes the proof of the second
equality.

The following is the main technical result of this section

\begin{theorem}\label{t2}
Let $t >5 $ and $f,\ g\in L^{\infty}({\overline{\DD}})$ are given. Suppose
moreover, that $g\in C^{\infty} (\overline{\DD})$ and that 
$$
\mbox{inf}\ \{ ||z-\xi || | \mid z\in \mbox{supp $f$ and $\xi \in$supp $g$}
\} > \epsilon 
$$
for some positive number $\epsilon$. Then both $T_f T_g$ and $T_g T_f$
are of trace class and Tr$ ([ T_f ,T_g ]) =0$.
\end{theorem}
Proof. We use $\partial$ to denote the unbounded operator
$$
{\mathcal H}_t  \ni h \rightarrow
\partial_z h \in {\mathcal H}_t .
$$
defined on the subspace of holomorphic functions $h$ such that their
first derivative is smooth up to the boundary of the disc and in
${\mathcal H}_t$, and by $\partial^{-1}$ the unique extension to a
bounded operator on ${\mathcal H}_t$ of
\begin{equation}\label{n}
z^n \rightarrow \frac{1}{n+1} z^{n+1} .
\end{equation}
It is easy to see that $\partial^{-1}$ is Hilbert-Schmidt, in fact,
since $||z^n ||_2 \sim O(1)$ as $n\rightarrow \infty$, the
characteristic values of $\partial^{-1}$ are of the order O($n^{-1}$).
Moreover $\Id - \partial \partial^{-1}$ is of finite rank. Since we can
write 
$$
T_f T_g |_{\C [z]}=T_f T_g \partial \partial^{-1} + 
                   T_f T_g (1-\partial \partial^{-1}),
$$
to prove that $T_f T_g$ is trace class it is sufficient to show that 
the densely defined operator $T_f T_g \partial$ has a (unique) extension
to a Hilbert-Schmidt operator on ${\mathcal H}_t$. Suppose first that
$h$ and $\partial h$ both belong to ${\mathcal H}_t$ and are smooth up
to the boundary of $\DD$. Given an $a \notin \mbox{ supp}(g)$, we have
$$
\begin{array}{lll}
 (T_g \partial h) (a)& = & \frac{t-1}{2\pi i} \int_{\DD} d\mu_t (b)
                     \frac{g(b)}{(1-a\bar{b})^t} \partial_b h (b) \\
                     & = &- \frac{t-1}{2\pi i} \int_{\DD}
          \frac{g(b)(1-|b|^2)^{t-2}}{(1-a\bar{b})^t}dh d\bar{b}  \\
                     &=& \frac{t-1}{2\pi i} \int_{\DD}
\partial_b \left( \frac{g(b)(1-|b|^2)^{t-2}}{(1-a\bar{b})^t} \right) h d\bar{b},
\end{array}
$$
where we used Stokes theorem and the fact that the integrand is smooth
and vanishes at the boundary of $\DD$. But this
implies that the densely
defined operator $T_f T_g \partial$ is in fact given by an integral
operator with kernel
$$
K(z,\xi)= const \int_{\DD \times \DD} d\lambda (a,b)
\frac{(1-|a|^2)^{t-2}(1-|b|^2)^{t-3}}{(1-\bar{a}z)^t (1-\bar{\xi}b)^t}
     \frac{F(a,b)}{(1-\bar{b}a)^t}
$$
where $d\lambda$ is the Lebesque measure on $\DD \times \DD$ and $F$ is
an $L^{\infty}$-function which vanishes on a neighbourhood of the
diagonal in $\overline{\DD} \times \overline{\DD}$ given by 
$\{ (a,b) \mid |a-b| >\epsilon \}$.
In particular,
$$
\sup_{a,b} |\frac{F(a,b)}{(1-\bar{b}a)^t}| <\infty
$$
and, by Cauchy-Schwartz inequality,
$$
|K(z,\xi )|^2 \leq const (Vol\ ( \DD , d\lambda )^2 
 \int_{\DD \times \DD}  d\lambda (a,b)
|\frac{(1-|a|^2)^{t-2}(1-|b|^2)^{t-3}}{(1-\bar{a}z)^t
(1-\bar{\xi}b)^t}|^2 .
$$
To estimate the $L^2$-norm of $K(z, \xi )$ we can first integrate over z
and $\xi$ which gives, in view of (\ref{5}), the estimate
$$
||K||_2^2 \leq const \int_{\DD \times \DD} (1-|a|^2)^{t-4}(1-|b|^2)^{t-6}
d\lambda (a,b)
$$
which is finite for $t>5$.

To finish the proof $T_g T_f = (T_{\bar{f}} T_{\bar{g}})^*$ and hence is
also trace class by applying the above argument to $T_{\bar{f}}
T_{\bar{g}}$. As a direct consequence we get $
Tr [T_f , T_g ]=0.
$

\section{$\Gamma$-invariant Toeplitz operators}

Let $\Gamma$ be a countable, icc and discrete subgroup of PSL(2,$\R$).

The von Neumann algebra $(\pi_t (\Gamma )^{''})$ is a $II_1$ factor with
unique normal normalized trace $\tau^{'}$ given by 
$$
\tau^{'} (\pi_t (\gamma ))=0
$$
for $\gamma \neq e$. Its commutant 
${{\mathcal A}_t}=(\pi_t (\Gamma ))^{'}$ is a factor of type II. We will
assume from now on that $\Gamma$ has infinite covolume in $\DD$, i. e. 
$$
M =\DD /\Gamma
$$
is an open Riemannian surface which can (and will) be thought of as an
open subset of an ambient closed Riemannian surface $N$. We assume
moreover that $M$ has finitely
many boundary components (the boundary in $N$), hence 
$\partial M = \cup_i C_i$, a finite union of disjoint
smooth simple closed contractible curves in $N$. In this case $\mathcal
M$ is a $II_{\infty}$ factor with a unique (up to the normalisation)
positive normal trace $\tau$. By general theory 
(see \cite{Ha}, \cite{C-1},\cite{EN}) there exists a
unique, normal, semifinite operator-valued weight
$$
E: {\mathcal B}({\mathcal H}_t )\rightarrow {\mathcal{M}}
$$
such that, for $A \in {\mathcal L}^1 (Tr)$ {\em in the domain of $E$},
$$
\tau \circ E (A) = Tr (A).
$$
$E$ is uniquely determined by the equality of normal linear functionals
$$
m \rightarrow \tau (E(A)m)=Tr(Am)
$$
for $A \in {\mathcal L}^1 (Tr)$ and $m\in {\mathcal A}_t$. 
Below we list some of the properties of $E$ used later.

A vector $\xi \in \HT$ is called $\Gamma $-bounded if the densely defined 
map
$$
l^2 (\Gamma )\ni \{ c_\gamma \}_{\gamma \in \Gamma} 
      \stackrel{R_{\xi}}{\rightarrow}\sum_{\gamma \in \Gamma} c_\gamma
\pi_t (\gamma^{-1} )\xi \in \HT
$$
is bounded. Let $p_\xi$ denote the orthogonal projection onto the
one dimensional subspace spanned by vector $\xi$. Then it is easy to see
that 
$$
R_\xi R_\xi^* = \sum_{\gamma \in \Gamma} \pi_t (\gamma ) P_{\xi} \pi_t
(\gamma^{-1} )
$$
i. e. $\xi$ is $\Gamma$-bounded precisely in the case when the sum
$\sum_{\gamma \in \Gamma} \pi_t (\gamma ) p_\xi \pi_t (\gamma^{-1} )$
converges in the
strong operator topology to a bounded oprator on $\HT$, in fact equal to
$E(p_\xi )$ and in this case $\tau (E(p_\xi ))=Tr (p_\xi )=1$. 

Let us introduce the following.

\begin{definition}
A bounded operator $A$ is called $\Gamma$-bounded if the sums
$$
\sum_{\gamma \in \Gamma}\pi_t (\gamma ) A \pi_t
(\gamma^{-1} )
$$
converge in the strong operator topology.
\end{definition}

Let $A$ be a positive trace class operator of the form
$$
A x = \sum_i \lambda_i p_{\xi_i} 
$$
where $\{ \xi_i \}$ is an orthonormal system in
$\HT$. $A$ is in the domain of $E$ if it is $\Gamma$-bounded and in this case 
$$
E(A)= \sum_{\gamma \in \Gamma} \pi_t (\gamma ) A \pi_t
(\gamma^{-1} ) \ \mbox{ and } \ \tau (E(A))=Tr A.
$$

\begin{proposition}\label{gcompact}
Let $f_0 $ be an $L^\infty$ function $\DD$ satisfying the conditions:
\begin{itemize}
\item the euclidean distance from the essential support of
$f_0$ to $\partial \DD$ is strictly positive;
\item the sum $\sum_{\gamma \in \Gamma} f_0 \circ \gamma$ is locally
finite. 
\end{itemize}
If, moreover, $t>2$, the associated 
$\Gamma$-invariant $L^\infty$-function $f= \sum_{\gamma \in \Gamma} f_0
\circ \gamma$ on $\DD$ satisfies 
$$ 
T_f \in {{\mathcal A}_t}\cap{\mathcal L}^1 (\tau)
$$ 
\end{proposition}
Proof. Since $f$ is $\Gamma$-invariant and in $L^\infty (\DD )$, $T_f
\in {{\mathcal A}_t}$. for the
rest of the claim it is sufficient to look at $f_0$ positive. 
But then $T_{f_0}$ is a positive operator with smooth kernel, hence, by
Lidskii theorem, it
is of 
trace class. By the second assumption it is $\Gamma$-bounded
and hence, according to the remarks above, it is in the
domain of $E$ and 
$$
E(T_{f_0}) = T_{\sum_{\gamma } f_0 \circ \gamma} =T_f .
$$
In particular
$$
\tau (T_{f})= \tau (E(T_{f_0}))=Tr( T_{f_0})<\infty 
$$
as claimed.

\bigskip

Let us introduce some notation connected with fundamental domains for
the action of $\Gamma$ on $\DD$. Suppose we choose points $P_i$ on
$\partial M$, one on each connected component $C_i$. We'll call this a
{\bf cut} of $M$. To each such cut we can associate a fundamental domain
$F$ such that the chosen points are in bijective correspondence with
end-points of the intervals $\overline{F}\cap \partial \DD$ 

From now on $F$ will (unless explicitly stated to the contrary) denote a generic
fundamental domain for $\Gamma$ on $\DD$.

\bigskip

Our goal is to compute the $\tau$-trace of commutators of the form $[T_f
,T_g]$, where $f$ and $g$ are sufficiently general $\Gamma$-invariant
functions on $\DD$. To see what is the problem, suppose first that
$f_0,g_0 \in C^{\infty} (F)$ satisfy supp$f_0
\subset
F^{\mbox{int}}$ and supp$g_0 \subset
F^{\mbox{int}}$. Let $f =\sum_{\gamma \in \Gamma} f_0 \circ \gamma$ and
$g=\sum_{\gamma \in \Gamma} g_0 \circ \gamma$ be the corresponding
$\Gamma$-invariant function on $\DD$. Looking at kernels, we obtain that
 (\cite{Ra}, \cite {NN}) 
$$
\tau (T_{f^2} -T_f^2 )=\frac{1}{2} (\frac{t-1}{2\pi i})^2 \int_{\DD \times
F} (f(a)-f(b))^2
\delta^t (a,b) d\mu_0
$$
and hence, by an application of the polarisation identity, 
$$
\tau (T_{fg}-T_{f} T_g )= (\frac{t-1}{2\pi i})^2  \int_{\DD \times F}
f(a)(g(a)-g(b))\delta^t (a,b) d\mu_0.
$$ Note that the right hand side is bounded by
$$
\frac{1}{2} (\frac{t-1}{2\pi i})^2 \int_{\DD \times
F} \vert a-b\vert^2\delta^t (a,b) d\mu_0,$$ which is convergent by \cite{Arazy}.

Hence $$\tau ([T_{f} ,T_{g} ] )=
(\frac{t-1}{2\pi i})^2  \int_{\DD \times F}
(f(a)g(b)-f(b)g(a))\delta^t (a,b) d\mu_0,
$$
Consequently, since the formula in Theorem 1 extends by 
continuity for functions $f,g$
that are smooth and $\Gamma$-invariant (by replacing ${\DD \times
\DD}$ by ${\DD \times
F})$, it follows that, for such $f$ and $g$,
$$
(\frac{t-1}{2\pi i})^2  \int_{\DD \times F}
(f(a)g(b)-f(b)g(a))\delta^t (a,b) d\mu_0 =
\frac{1}{2\pi i} \int_{\DD}
d(f_0 )dg =\frac{1}{2\pi i} \int_F dfdg .
$$
But it is not obvious from the outset neither that $\tau ([T_{f} ,T_{g} ] )$
is in the domain of $\tau$ nor that its trace 
is approximated by the trace of commutators of
Toeplitz operators associated to functions of the form $f =\sum_{\gamma
\in \Gamma} f_0 \circ \gamma$ and $g=\sum_{\gamma \in \Gamma} g_0 \circ
\gamma$ with $f_0$ and $g_0$ supported away from the boundary!

\begin{theorem}\label{t3}
Let $f$ and $g$
 be two smooth functions on $\overline{M}$ (i.e. continuous with all
their derivatives up to the boundary of $M$). We will use the same
notation to denote their representatives as $\Gamma$-invariant functions
on $\DD$. Suppose that $t>5$. Then $[T_f ,T_g ]$ is in ${\mathcal
M}\cap{\mathcal L}^1
(\tau)$ and
$$
\tau ([T_f ,T_g ] )=\frac{1}{2\pi i} \int_F df dg .
$$
\end{theorem}
Proof.

Let us begin with the following observations.
\begin{enumerate}
\item
Suppose that $h_0 ,\ldots h_n$ is a finite family of functions on $\DD$
satisfying the conditions of the proposition \ref{gcompact} and we set $
A= T_{h_0} \ldots T_{h_n} $. Since
$$
|A^*|^2 \leq \left( \prod_{i\neq 0} ||h_i ||_\infty  \right) T_{|h_0 |^2},
$$ 
the averages 
$\sum_{\gamma} \pi_t (\gamma ) |A^*|^2 \pi_t (\gamma^{-1})$ converge
in the strong
operator topology to $E(|A^*|^2)$. Moreover, for any normal linear
functional
$\psi$ on ${\mathcal B }(\HT )$, 
$$
\sum_\gamma \psi (\pi_t (\gamma ) |A^*|^2 \pi_t
(\gamma^{-1}))= \psi (E(|A^*|^2)).
$$
To see the equality it is sufficient to consider positive $\psi$, but
then all that is involved is an exchange of the order of summation for a
double series consisting of positive terms.
\item
There exists a smooth partition of unity of
$\DD$
of the form $\sum_i \phi_i$ where $\phi_i$ are smooth, positive
functions such that, for each $i$, the family of functions $\{ \phi_i \circ
\gamma \}_{\gamma \in \Gamma}$ is localy finite. To see this it is
sufficient to notice that, for any disc $\DD_\epsilon = \{ z\in \DD | |z|
\leq \epsilon \}$ with $\epsilon <1$, the number of $\gamma \in \Gamma$
such that $\gamma (\DD_\epsilon )\cap \DD_\epsilon \neq \emptyset$ is
finite. This follows from the fact that $\Gamma$ is discrete as a
subgroup of PSL(2,$\R$.

This implies that the
set of normal linear functionals $m \rightarrow Tr( BmA^*)$ with $A$ and
$B$ as above is
total in $\mathcal{B} (\HT )_*$.
\end{enumerate}

The proof of the theorem will be done in two steps.

\medskip

\noindent {\em First part.}

\medskip

Suppose that we are given a fundamental domain $F$ for $\Gamma$
such that 
$$
f=\sum_\gamma f_0 \circ \gamma ,\ g=\sum_\gamma g_0 \circ \gamma
$$
where $f_0$ and $g_0$ are both smooth on $\overline{\DD}$ and their
supports have positive Euclidean distance to the complement of $F$ in
$\DD$.

We will compute
$$
Tr (A^* [T_f ,T_g ] A )
$$
where $A= T_{h_0} \ldots T_{h_n}$.
The operator under trace has a smooth kernel and the integral of its
restriction to the diagonal has (up to a constant) the form
$$
\int_\DD d\mu_t (d)\int_\DD d\mu_t (c)\int_\DD d\mu_t (b)\int_\DD d\mu_t (a)
\frac{A(a,d)(f(b)g(c)-g(b)f(c))}{(1-\bar{a}b)^t
(1-\bar{b}c)^t (1-\bar{c}d)^t (1-\bar{d}a)^t}
$$
where $A(a,d)$ is a smooth kernel with support of strictly positive
euclidean distance from $\partial \DD \times \DD \cup \DD \times
\partial \DD$.
Hence the function
$$
F(a,b,c,d) =(1-|b|^2)^{t/2} (1-|c|^2)^{t/2} \frac{A(a,d)(f(b)g(c)-g(b)f(c))}{(1-\bar{a}b)^t
(1-\bar{b}c)^t (1-\bar{c}d)^t (1-\bar{d}a)^t} 
$$
is uniformly bounded on $\DD^4$ (the only singularity in the
denominator appears for $b=c$ and it is controlled by the fact that
$|\delta|\leq 1$) and our
integral can be written as the integral of $L^{\infty}$-function
$F(a,b,c,d)$ with respect to the finite measure
$$
d\omega = d\mu_t \otimes  d\mu_{t/2} \otimes d\mu_{t/2}
\otimes d\mu_t.
$$
Hence
$$
\int F d\omega = \sum_{\gamma \in \Gamma} \int_{\DD \times \gamma (F)
\times \DD \times \DD } F d\omega ,
$$
i. e.
$$
Tr (A^* [T_f ,T_g ]A )=\sum_{\gamma} Tr (A^* (T_{f_0
\circ \gamma } T_{g}-T_{g_0 \circ \gamma }T_f
)A ).
$$

Since $T_{f_0} T_{g}-T_{g_0}T_f$ is of trace class, we can exchange the
summation over $\gamma \in \Gamma$ with the trace and get the identity
$$
Tr (A^* [T_{f} ,T_{g} ]A)= Tr \left( \sum_{\gamma}\pi_t (\gamma )A A^*
\pi_t (\gamma )^{-1}(T_{f_0} T_{g}-T_{g_0}T_f
)\right) .
$$

But this shows that
$$
\tau (E(|A^*|^2)[T_f ,T_g ])= 
            \tau (E(|A^*|^2)E( T_{f_0} T_g -T_{g_0} T_f)).
$$
Since for any $m\in {\mathcal A}_t$
$$
\tau (E(A A^* )m)=Tr((AA^* )m)
$$
and the set of such linear functionals is separating for $\mathcal B (\HT )$, 
we get
$$
E( T_{f_0} T_g -T_{g_0}
T_f) = E( T_{f} T_{g_0} -T_{g} T_{f_0}) =[T_f ,T_g ].
$$
But, since
$$
T_{f} T_{g_0} -T_{g} T_{f_0}=[ T_{f_0}, T_{g_0}] +T_{f-f_0}
T_{g_0} -T_{g-g_0} T_{f_0}
$$ 
is of trace class, $[T_f ,T_g ]$ is in ${\mathcal
M}\cap{\mathcal L}^1 (\tau)$ and 
$$
\tau ([T_f,T_g])= \frac{1}{2} Tr ([T_{f_0}, T_g ]
-[T_{g_0}, T_f ])
$$

By the theorem 3
$$
\tau ([T_f ,T_g ] )=\frac{1}{2\pi i} \int_{\DD}
d(f_0 )dg =\frac{1}{2\pi i} \int_F dfdg .
$$

\medskip

\noindent {\em Second part.}

\medskip

By the proposition 1 we can assume that
both $f$ and $g$ are, as functions on $\overline{M}$, supported on a
neighbourhood
of the boundary of $\overline{M}$ diffeomorphic to 
$( \cup_i  C_i )\times ]-\epsilon ,0]$. Now, using partition of unity, we
can split both $f$ and $g$ into finite sums
$$
f=\sum_k f_k ,\ g=\sum_s g_s
$$
so that for any pair of indices $(k,s)$ there are open intervals $I_{k,s}^i$ of
non-zero length on each of the boundary components $C_i$ such that both $f_k$
and $g_s$ vanish on $(\cup_i I_{k,s}^i \times ])-\epsilon ,0]$ -
possibly with a smaller, but still positive value of $\epsilon$. But
then, choosing a cut of $M$ given by a choice of points $P_i$ in the interior
of $I_{k,s}^i$ will provide us with a fundamental domain $F_{k,s}$ such
that the conditions of the first part of this proof hold for $(f_s ,g_k
,F_{k,s} )$ and hence $[T_{f_k},T_{g_s}]\in {\mathcal L}^1 (\tau )$ and
$$
\tau ([T_{f_k} ,T_{g_s} ] ) =\frac{1}{2\pi i} \int_{F_{k,s}} dfdg
$$ 
To complete the proof note that the expression
$dfdg$ for $\Gamma$-invariant functions is $\Gamma$-invariant,
hence the integral $\int_{F} dfdg$ is independent on the choice
of the fundamental domain and the result follows.

\begin{corollary}
Suppose that $f$ and $g$ are smooth functions on $\overline{M}$. Then
$$
\tau ([T_{f} ,T_{g} ] )=\int_{\partial M} fdg.
$$
\end{corollary}
Proof. This follows immediately from the fact that under the
natural diffeomorphism $F \setminus \partial F$ the integral $\int_{F}
dfdg$ becomes identified with $\int_M dfdg$ and the Stokes theorem.

\bigskip

\noindent {\bf Remarks}
\begin{enumerate}
\item
Virtually the same proof shows that, for $f$ and $g$ smooth on
$\overline{M}$, the operator $T_f T_g - T_{fg}$ is in ${\mathcal
M}\cap{\mathcal L}^1
(\tau)$ and 
$$
\tau (T_f T_g - T_{fg} )=(\frac{t-1}{2\pi i})^2  \int_{\DD \times F}
f(a)(g(b)-g(a))\delta^t (a,b) d\mu_0 (a,b).
$$
\item 
All of the results above can be easily extended to the case when $f$ and
$g$ are in $L^infty (M)$ and Lipschitz with exponent one in a
tubular neighbourhood of $\partial M$ in $\overline{M}$.
\end{enumerate}

\section{$\Gamma$-Fredholm operators}

Let ${\mathcal T}_{\Gamma}$ denote the C*-subalgebra of ${\mathcal A}_t$
generated by Toeplitz operators $T_f$ with Toepltz symbol $f\in
C(\overline{M})$ and denote by ${\mathcal K}_{\Gamma}$ the C*-ideal
generated by elements in ${\mathcal L}^1 (\tau) \cap {\mathcal
T}_{\Gamma}$. An element $A$ of ${\mathcal A}_t$ is called $\Gamma$-Fredholm
if it has an inverse, say $R$, modulo ${\mathcal K}_{\Gamma}$ and, 
in this case, the commutator $[A,R]$ has well-defined trace 
$$
\mbox{$\Gamma$-index of $A$}=\tau ([A,R]) 
$$ 
which depends only on the class of $A$ in 
$K_1 ( {\mathcal T}_{\Gamma}/{\mathcal K}_{\Gamma} )$.

\begin{remark}
The number'' $\Gamma$-index of $A$'' is also known as 
{\em Brauer index} of $A$.
\end{remark}
According to the proposition 1, there exists a surjective
continuous map

\begin{equation}\label{n}
q: C(\partial M) \rightarrow {\mathcal T}_{\Gamma}/{\mathcal K}_{\Gamma}
\end{equation}
sending function $f|_{\partial M}$ to $T_f$ mod ${\mathcal K}_{\Gamma}$
- this map is well defined since $|| T_f || =||f||_{\infty}$.

\begin{theorem} Let $\Gamma$ be a countable, discrete, icc subgroup of
PSL(2,$\R $) such that $\DD /\Gamma$ has infinite covolume and $M=\DD
/\Gamma$ is an open Riemannian surface with finitely many boundary
components. Assume that $t>5$ and that the trace $\tau$ on ${\mathcal A}_t$
is normalized by its value on a Toeplitz operators
$T_f$ with symbols $f\in C_c^{\infty} (M)$ by
$$
\tau (T_f) = \frac{t-1}{2\pi i} \int_F f(z) d\mu_0 (z).
$$
The following holds.

\noindent 1. For any function $f\in C(\overline{M})$ such that $f|_{\partial
M}$ is invertible, the operator $T_f$ is $\Gamma$-Fredholm
and its $\Gamma$-index is equal to the sum of the winding numbers of
restriction of $f$ to the boundary of $M$.

\noindent 2. The map $q: C(\partial M) \rightarrow {\mathcal
T}_{\Gamma}/{\mathcal K}_{\Gamma}$ is injective and yields a nontrivial
extension 
$$
0\rightarrow {\mathcal K}_{\Gamma} \rightarrow {\mathcal
T}_{\Gamma} \rightarrow C(\partial M) \rightarrow 0.
$$

\end{theorem}
Proof. 

\noindent {\em Step 1.}

Suppose that $f\in C^\infty (\overline{M})$ be invertible on the
boundary of $M$. Then, for any function g smooth in the closure of $M$, 
and such that supp$(1-fg)\subset M$ the theorem 3 gives $1-T_f T_g \in
{\mathcal K}_{\Gamma} $ and hence $T_f$ is $\Gamma$-Fredholm.

Let now $P$ and $Q$ be two non-commutative polynomials in $(\bar{z},z)$.
Then, again by the theorem 3 and its corollary,
$$
\tau ([P({T_f}^* ,T_f ),Q({T_f}^* ,T_f )])=
\tau([ T_{P(\bar{f} ,f)}, T_{Q(\bar{f},f)}])=\frac{1}{2\pi i}\int_{\partial
M} P(\bar{f} ,f)dQ(\bar{f} ,f).
$$
On the other hand, by Carey-Pincus formula for traces of commutators
(see \cite{CP1}),
$$
\tau ([P({T_f}^* ,T_f ),Q({T_f}^* ,T_f )])= \int_{|z|<||T_f ||} \{
P,Q \} d\nu
$$
where $d\nu$ is a finite measure supported on the convex hull of the
essential spectrum of $T_f $mod(${\mathcal K}_{\Gamma}$) and, since
$T_f$ is $\Gamma$-Fredholm, there exists an open ball $B_\epsilon$
around the origin such that 
$d\nu|_{B_\epsilon} = cd\lambda$, where $ 2\pi i c = \Gamma
\mbox{-index of $T_f$}$.

If we set $d\nu = d\nu|_{B_\epsilon} +d\nu_1$, the two formulas above give 
$$
\frac{1}{2\pi i}\int_{\partial
M} P(\bar{f} ,f)dQ(\bar{f} ,f)=\frac{1}{2\pi i}\Gamma
\mbox{-index of $T_f$}
\int_{|z|\leq \epsilon}dPdQ +\int_{\epsilon \leq |z| \leq ||T_f ||}\{
P,Q\} d\nu_1 .
$$
Applying Stokes theorem, we get the equality
$$
\frac{1}{2\pi i}\int_{\partial
M} P(\bar{f} ,f)dQ(\bar{f} ,f) = \frac{1}{2\pi i}(\Gamma
\mbox{-index of $T_f$})
\int_{|z|= \epsilon}PdQ +\int_{\epsilon \leq |z| \leq ||T_f ||}\{
P,Q\} d\nu_1 .
$$
If we now set $P(\bar{z},z)=z$ and approximate 
$\frac{1}{z}$ uniformly on the annulus $\epsilon \leq |z| \leq ||T_f ||$
by polynomials $Q$, since both sides are continuous in the uniform
topology on $C( \epsilon \leq |z| \leq ||T_f || )$ we get, in the limit,
$$
\frac{1}{2\pi i}\int_{\partial
M} f^{-1} df = \frac{1}{2\pi i}(\Gamma
\mbox{-index of $T_f$})
\int_{|z|= \epsilon}z^{-1}dz =( \Gamma
\mbox{-index of $T_f$}).
$$

\noindent {\em Step 2.}

We will now prove injectivity of $q$. To this end it is enough to show
that, for any open interval $I\subset \partial M$, we can find a function
in $C(\overline{M})$ which is zero when restricted to $\partial M
\setminus I$ but for which the corresponding toeplitz operator $T_f$ is
not in ${\mathcal K}_{\Gamma}$. But given such an interval, we can
easily find a smooth function $f$ such that $f|_{\partial M}$ is
supported within $I$ and such that the winding number of $1+f$ on
$\partial M$ is nonzero. But then, by step 1 above, the $\Gamma$-index
of $1+T_f$ is nonzero and hence $T_f $ is not an element of $
 {\mathcal K}_{\Gamma}$,
which proves the second part of the theorem.

Now, let $f$ be continuous and with invertible 
restriction to the boundary of $M$. By part two
of the theorem, this implies that the image of
$T_f$ in ${\mathcal T}_{\Gamma} /{\mathcal K}_{\Gamma}$ is $f|_{\partial
M}$ and hence invertible, i. e. $T_f$ is
$\Gamma$-Fredholm, and the
formula for its $\Gamma$-index follows from the fact that it is a
functional on $K_1 ({\mathcal T}_{\Gamma} /{\mathcal K}_{\Gamma})$ and
hence homotopy invariant of the class of $f|_{\partial M}$ in $K_1
(C(\partial M))$.

The normalisation statement follows immediately from the proposition 1.

\end{document}